\documentclass[letterpaper,11pt]{article}
\usepackage[latin1]{inputenc}
\usepackage[T1]{fontenc}
\usepackage[spanish]{babel} 
\usepackage[pdftex]{hyperref}
\usepackage{amsmath,amsfonts,amssymb}
\usepackage{amsthm}
\usepackage{multicol} 
\usepackage{graphicx}
\usepackage{fancyhdr}
\usepackage[all]{xy}
\usepackage{pgf, tikz}
\usepackage{multicol}

\theoremstyle{plain}
\newtheorem{theorem}{Teorema}[section]{\bfseries}{\scshape}  
{\bfseries}{\itshape} 
\newtheorem{propo}{Proposici\'on}[section]{\bfseries}{\itshape} 
{\bfseries}{\itshape}
\newtheorem{defi}{Definici\'on}[section]{\bfseries}{\itshape} 

\theoremstyle{definition}
{\bfseries}{\rmfamily}
{\bfseries}{\rmfamily}
\newtheorem{exam}{Ejemplo}{\bfseries}{\rmfamily}

{\bfseries}{\rmfamily}
\newtheorem*{remark}{Observaci\'on}{\bfseries}{\rmfamily}
\newtheorem*{remarks}{Observaciones}{\bfseries}{\rmfamily}

\newcommand{\R}{\mathbb{R}}

\usepackage[left = 2cm, right = 2cm, top = 2cm, bottom = 2cm]{geometry}
\usepackage[margin=3cm, font=footnotesize, labelfont=bf]{caption}

\author{\sc   Yoceman Sifontes\\
	Departamento de Matemática, UCA, San Salvador, El Salvador\\
	Teléfono (503)  76042269\\
	E-mail: ysifontes@uca.edu.sv\vspace*{1cm}\\
		Dimas  Tejada\\
		Facultad de Ciencias Naturales y Matemática \\
		Universidad de El Salvador, El Salvador \\
		\textit{E-mail:} {\tt dimas.tejada@ues.edu.sv } }
\title{\bf Aproximación Dinámica al Problema de Hermite}
\date{Mayo de 2023}

\hypersetup{colorlinks=true,linkcolor=black}

\begin{document}

\maketitle

  \centerline{\textit{\bf Resumen}}  En \cite{GCF} se prueba que cualquier número irracional cuadrático tiene una representación como fracción continua,  infinita y periódica. En 1848, Charles Hermite mediante una carta Jacobi \cite{Per}, se preguntó, si este hecho se podría generalizar para números irracionales cúbicos. En esta investigación se aborda el problema de Hermite, mediante el estudio de funciones continuas por pedazos asociadas a Grupos Fuchsianos, particularmente una función $f_\Gamma$ asociada al Grupo Modular $\Gamma$, que fu\'e planteada en 1978 por Robert Edward y Caroline Series, con esta función se encuentra una representación como fracción continua infinita periódica para números irracionales cuadráticos e infinita no periódica para irracionales no cuadráticos. \\
  
\textit{Palabras Clave} - Aproximación Dinámica, Grupos Fuchsianos, Problema de Hermite.\\

 \centerline{\textit{\bf Abstract}}
In \cite{GCF} it is proved that any quadratic irrational number has a representation as a continuous, infinite and periodic fraction. In 1848, Charles Hermite through a letter
Jacobi \cite{Per} wondered if this fact could be generalized to cubic irrational numbers. In This research addresses the Hermite problem, through the study of continuous functions by pieces associated to Fuchsian Groups, particularly a function $f_\Gamma$ associated to the Group
Modular $\Gamma$, which was raised in 1978 by Robert Edward and Caroline Series, with this function it is find a representation as infinite periodic continued fraction for irrational numbers quadratics and non-periodic infinite for non-quadratic irrationals.\\


\vspace*{1cm}


\begin{center}\bf Introducci\'on\end{center} 
En \cite{BC} se establece: dado un grupo Fuchsiano $\Gamma$ finitamente generado que actúa en el disco de Poincaré $\mathbb{D}$, los autores describen una función continua por pedazos  definida en la frontera ideal $\partial \mathbb{D}=\mathbb{S}^1$, que refleja las propiedades dinámicas de la acción del grupo. Utilizando elementos de dinámica simbólica, se estudia la función $f_\Gamma$, con $\Gamma$ el Grupo Modular, asociando un espacio de s\'imbolos. La funci\'on $f_\Gamma$  se define por primera vez en \cite{BC} como una función continua por pedazos; para dicha función se establecen resultados referentes a la naturaleza de sus órbitas periódicas bajo $f_\Gamma$.  La representación simbólica del itinerario de un número real $x$ bajo la función $f_\Gamma$ permite, para n\'umeros irracionales cuadr\'aticos, asociar una representaci\'on como fracci\'on continua peri\'odica; mientras que, para n\'umeros irracionales no cuadr\'aticos, permite encontrar una representaci\'on como fracci\'on continua, infinita no peri\'odica bajo $f_\Gamma$.

\section{Preliminares}

\subsection{Fracciones continuas}
El término de fracción continua fue usado por primera vez en 1695 por J. Wallis, pero el uso de las fracciones continuas data desde los tiempos de Euclides. Para más detalles históricos, consultar \cite{DF}.
\begin{defi}\label{fcontinua}
	La representaci\'on de un n\'umero real $x$ como fracci\'on continua es 
	\[x=[x_0;x_1:x_2:x_3:\cdots ]=x_0+\cfrac{1}{x_1+\cfrac{1}{x_2+\cfrac{1}{x_3+\cdots } }},\]
	donde $x_0\in \mathbb{Z}$ y $\{x_i\}^\infty_{i=1}$,  es alguna secuencia de n\'umeros enteros no nulos.  Los números $x_0,\cdots, x_n$ son llamados elementos de la fracción continua. La representaci\'on como fracci\'on continua, puede ser finita o infinita.
	Si los $x_1,x_2,\cdots $ son enteros positivos, la fracción continua es llamada regular.
\end{defi}
\begin{defi}
	Una fracción continua finita se dice que es par o impar si tiene un número par o impar de elementos respectivamente.
\end{defi}
\begin{defi}
	Una fracción continua con una secuencia infinita de elementos $x_1,x_2,x_3,\cdots, $ es  denotada por $[x_0;x_1:x_2\cdots :x_k]$ si el límite \[ \lim_{k\rightarrow \infty}[x_0;x_1:x_2:\cdots :x_k]\] existe. 
\end{defi}

\begin{remark}
	La definici\'on puede variar \ref{fcontinua} dependiendo del contexto o intereses de trabajo, para tener una perspectiva m\'as amplica sobre la teor\'ia de fracciones continuas, consultar \cite{GCF}. 
\end{remark}
Dado un número real $x$, el siguiente teorema garantiza la existencia y unicidad para fracciones continuas  regulares, \cite{DF}.
\begin{theorem}
	\begin{enumerate}
		\item Para todo número racional, existe una única fracción continua par y una única representación continua impar.
		\item Para todo número irracional $x$ existe una única fracción continua regular, la cual converge a $x$ cuando $k\rightarrow \infty$.
	\end{enumerate}
\end{theorem}

\begin{exam}
	Para $x=\frac{10}{7}$, se tiene \[\frac{10}{7}=1+\cfrac{1}{2+\frac{1}{3}}=[1;2:3].\]
\end{exam}

\begin{exam}
	Para $\pi=3+\cfrac{1}{7+\frac{1}{15+\frac{1}{1+\cdots }}}=[3;7:15:1:292:1:1:1:2\cdots ]$
\end{exam}
En general, las fracciones continuas ofrecen una manera de conocer la irracionalidad de un número, ya que si  su desarrollo es infinito, entonces el número es irracional. Esta técnica fue utilizada por Euler, que determinó la fracción continua del número $e$, por otro lado, Johann Heinrich Lambert demostró por primera vez en 1761 la irracionalidad del número $\pi$ mediante el uso de fracciones continuas llamadas fracciones continuas generalizadas.

\subsection{Definición de la función $f_\Gamma$}

El estudio  de aplicaciones continuas por pedazos es un \'area relativamente nueva en matemática y existen  algunos resultados generales para cierto tipo  de transformaciones en diversos espacios. En este trabajo consideramos como espacio a $X=\R$ particionado en los conjuntos $(-\infty, -1)$, $(-1,0), (0,1)$ y $(1,+\infty)$.  La aplicación aparece en los  trabajos de Bowen-Series \cite{BC};  para estudiar su din\'amica, se define un espacio de s\'imbolos $\sum_{f_\Gamma}$ sobre la extensión de $f_\Gamma$ a $\widehat{\R}=\R\cup \infty$. El espacio de itinerarios posibles bajo $f_\Gamma$ permite identificar la existencia de \'orbitas peri\'odicas. Las aplicaciones a estudiar, en general, están asociadas a Grupos Fuchsianos, que son ciertos grupos de isometrías que preservan la orientación en el semiplano superior. Una excelente referencia para estudiar Grupos Fuchsianos es \cite{SK}.

\begin{defi}
	El Grupo Modular $\Gamma=\mbox{PSL}(2,\mathbb{Z})=\left\{\left. \dfrac{az+b}{cz+d}\ \right| a,b,c,d\in \mathbb{Z}, ad-bc=1\right\}$.
\end{defi}
El Grupo Modular es un grupo de gran importancia por sus conexiones con la Teor\'ia de N\'umeros, Topolog\'ia, Geometr\'ia Hiperb\'olica entre otras ramas de la matem\'atica.

\begin{defi}\label{funcionModular}
	Para $\Gamma=\langle x+1,-\frac{1}{x}\rangle $ el Grupo Modular, se define \[f_\Gamma(x)=\left\{
	\begin{array}{cc}
	x+1, & x\leq -1  \\
	-\dfrac{1}{x},& -1<x<1\\
	x-1,& x\geq 1\\
	\infty,&\infty.
	\end{array}\right.\]
\end{defi}

En esta sección, las sigiuentes definiciones suponen que $f: X\rightarrow X$ una función y $X$ tiene estructura de Espacio Topológico.

\begin{defi}
	A la colecci\'on de puntos $x, f(x), f^2(x),\cdots $ es llamada la \'orbita de $x$ y se denota por $\mathcal{O}^+(x)$.
\end{defi}
Si $f$ es un homeomorfismo, se define la \'orbita completa de $x$ como $\mathcal{O}(x)$, como la colecci\'on de puntos $f^n(x)$ con $n\in \mathbb{Z}$  y  a la \'orbita para los negativos se le denotar\'a por $\mathcal{O}^-(x)$.  

\begin{defi}
	Se dice $x$ es un punto fijo de $f$ si $f(x)=x$.   
	El punto $x$ es un punto peri\'odico de per\'iodo $n$ si $f^n(x)=x.$ Al menor entero positivo $n$ para el cual $f^n(x)=x$ le  llamaremos per\'iodo de $x$.  
\end{defi}

Denotaremos  por $\mbox{Per}_n(f)$ al conjunto de todos los  puntos peri\'odicos  de peri\'odo $n$  y por $\mbox{Fix}(f)$ al conjunto de  puntos fijos.  El conjunto de todas las iteraciones de un punto peri\'odico forma una \'orbita peri\'odica.\\

\begin{defi}
	Un punto $x$ es eventualmente peri\'odico  de periodo $n$, si $x$ no es peri\'odico pero existe $m>0$ tal que $f^{n+i}(x)=f^i(x)$ para todo $i>m.$ Esto es, $f^i(x)$ es peri\'odica para $i\geq m.$
\end{defi}

\subsection{Elementos de din\'amica simb\'olica}

Para m\'as detalles sobre din\'amica simb\'olica, ver \cite{LM} y \cite{RD}.

\begin{defi}
	Se denotar\'a por $\mathcal{A}$ a un conjunto finito o numerable de s\'imbolos que llamaremos {\bf alfabeto.} Usualmente $\mathcal{A}=\{0,\ldots, N-1\}$ \'o $\mathcal{A}=\{0,1,2\ldots \}$.
\end{defi}

\begin{defi}
	Si $\mathcal{A}$ es un alfabeto finito, entonces el  corrimiento completo es la colecci\'on de todas las secuencias bi-infinitas de s\'imbolos de $\mathcal{A}$. El  $r-corrimiento$ es el corrimiento completo sobre el alfabeto $\{0,1,\ldots, r-1\}.$
	
	Elcorrimiento completo es denotado por \[\mathcal{A}^\mathbb{Z}=\{ x=(x_i)_{i\in\mathbb{Z}}\mid x_i\in \mathcal{A},\forall i\in \mathbb{Z}\}=\sum.\]
\end{defi}

\begin{defi}
	La aplicaci\'on {\bf corrimiento unilateral} \[\sigma: \sum\rightarrow \sum \]
	se define 
	\[\sigma(x_i)=(x_{i+1}),\quad \forall i\in \mathbb{Z}\quad \mbox{o}\quad  \,\mathbb{N}.\]
\end{defi}

\begin{defi}
	Un espacio de corrimiento o simplemente corrimiento es un subconjunto  $X$ de  $\mathcal{A}^\mathbb{Z}$ tal que $X=X_\mathcal{F}$ para alguna colecci\'on $\mathcal{F}$ de bloques prohibidos sobre $\mathcal{A}$.
\end{defi}

\begin{defi}
	Un corrimiento de tipo finito es un corrimiento que puede ser descrito por un conjunto finito de bloques prohibidos, equivalentemente un corrimiento  $X$ que tiene la forma $X_\mathcal{F}$ para alg\'un conjunto finito $\mathcal{F}$ de bloques.
\end{defi}

\begin{exam}
	Por ejemplo, si $f:X\rightarrow X$ y $\mathcal{P}=\{P_0,P_1,P_2,P_3,P_4\}$.

	Si $f: X\rightarrow X$ es una funci\'on y $\mathcal{P}=\{P_i\}^{4}_{i=0}$ es una colecci\'on de subconjuntos de $X$ con interiores disjuntos y $X=\bigcup P_i$. 
	
	\begin{itemize}
		\item El conjunto $\mathcal{A}=\{0,1,2,3, 4\}$ es el alfabeto.
		\item $\mathcal{A}^\mathbb{N}=\{0,1,2,\ldots, 4\}^\mathbb{N}=\left\{(s_i)^\infty_{i=0}\mid s_i\in \mathcal{A}\right\}=\sum$.
		\item Si $x\in X$, $I(x)$ es la secuencia de $0's,\ldots 4's$ donde $s_j=0$ si $f^j(x)\in P_0$,\ldots, $s_j=4$ si $f^j(x)\in P_{4}$.
		\item $\sum_f$ la colecci\'on de los $I(x)$ permitidos bajo $f$.
		\item El corrimiento unilateral:
		\begin{align*}
		\sigma:\sum_f&\longmapsto\sum_f\\
		(s_0s_1s_2s_3\ldots )&\longmapsto(s_1s_2s_3\ldots )
		\end{align*}
	\end{itemize}

	\begin{figure}[h!]
		\centering
		\begin{center} 
			\includegraphics[width=10cm]{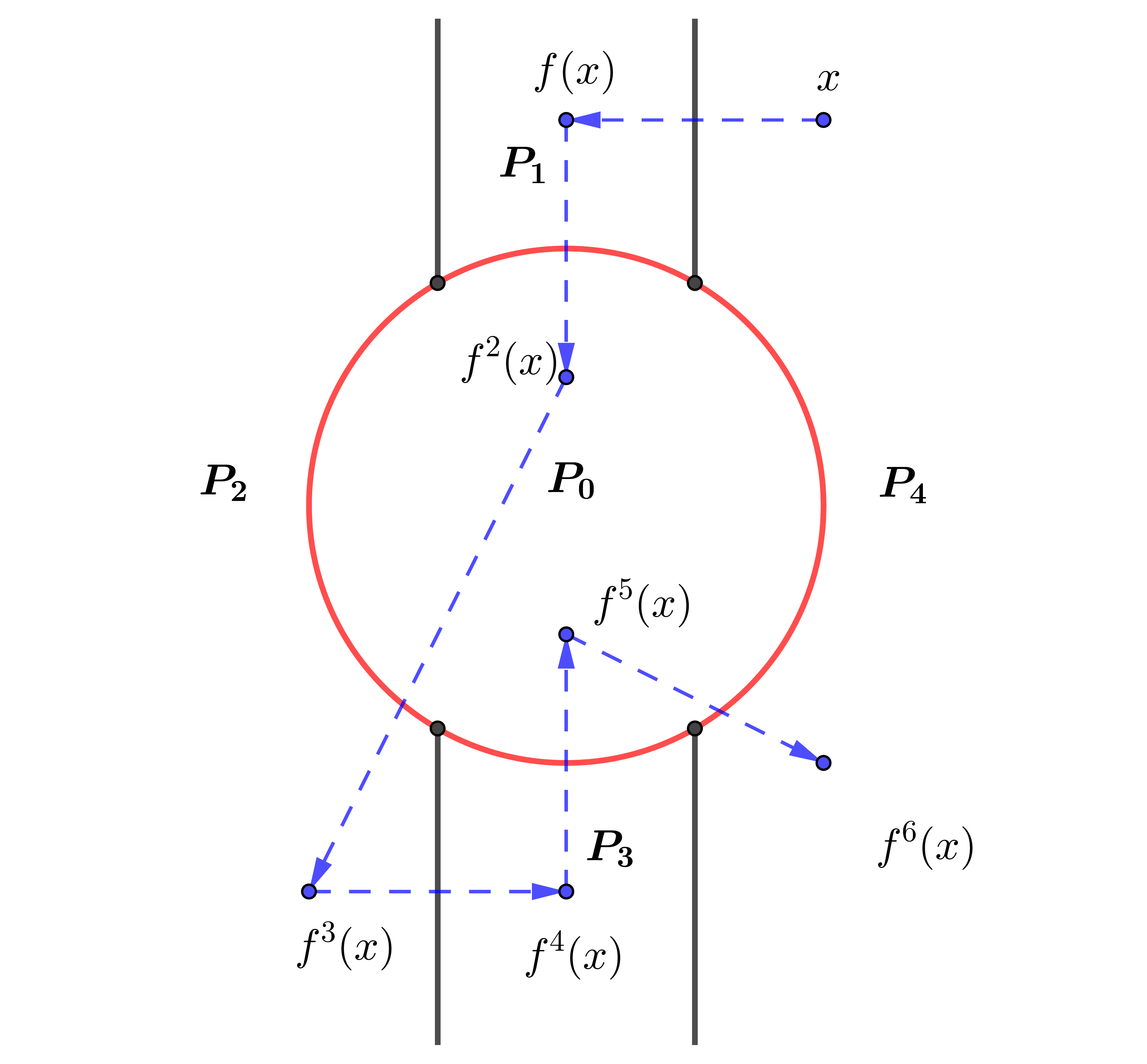}
			\caption{El itinerario para $x$ es $(4102304\ldots )$}
		\end{center}
		
	\end{figure}
\end{exam}

En general, para un espacio $\mathcal{A}^\mathbb{N}$ para alg\'un alfabeto  $\{0,1,2,\ldots, N-1 \}$ se tiene:

\begin{propo}
	El espacio $\sum$ es un espacio m\'etrico,  con la m\'etrica 
	\[d[s,t]=\sum^\infty_{i=0}\dfrac{|s_i-t_i|}{N^i}.\]
\end{propo}

\begin{propo}\label{proposicion221}
	Sean $s$ y $t$ elementos de $\sum$ y sup\'ongase que $s_i=t_i$ para $i=0,1,2,\ldots, N-1$. Entonces \[d[s,t]\leq \dfrac{1}{N^{n}}\]
\end{propo} 

\begin{proof} 
	Se tiene que si $s_i=t_i$ para $i\leq n$
	\begin{align*}
	d[s,t]&=\sum^\infty_{i=0}\dfrac{|s_i-t_i|}{N^i}\\
	&=\sum^n_{i=0}\dfrac{|s_i-t_i|}{N^i}+\sum^\infty_{i=n+1}\dfrac{|s_i-t_i|}{N^i}\\
	&\leq \sum^\infty_{i=n+1}\dfrac{N-1}{N^i}=\dfrac{N-1}{N^n}\sum^\infty_{i=1}\dfrac{1}{N^n}=\dfrac{N-1}{N^i}\cdot \dfrac{1}{N-1}=\dfrac{1}{N^n}.
	\end{align*}
\end{proof}

En general,  $\mathcal{A}^\mathbb{N}$ es un conjunto de Cantor, es decir, un conjunto compacto, totalmente disconexo sin puntos aislados.

\begin{propo}\label{Itinerario}
	La aplicaci\'on de itinerario $I: X\rightarrow \sum$ es continua.
\end{propo}

\begin{proof}
	Sea $f: X\rightarrow X$ una funci\'on continua y $\mathcal{A}=\{0,1,\cdots, N-1\}$ alg\'un alfabeto sobre una partici\'on $\{P_i\}^N_{i=1}$, para $x\in X$ sup\'ongase que $I(x)=s_0s_1s_2\ldots$, mostraremos que $I$ es continua en $x$.\\
	Sea $\epsilon>0$, se elige $n$ de tal forma que $\dfrac{1}{N^n}<\epsilon$ y se consideran los conjuntos 
	
	\[P_{t_0t_1\ldots t_n}=P_{t_0}\cap f^{-1}(P_{t_1})\cap \ldots f^{-n}(P_{t_n}),\] para todas las combinaciones de $t_0t_1\ldots t_n$, estos subconjuntos son disjuntos y su uni\'on cubre a $X$. Existen $N^{n+1}$  de tales subconjuntos  y $P_{t_0t_1\ldots t_n}$ es uno de ellos. Luego, se elige $\delta$ de tal forma que $|x-y|<\delta $, donde $y\in X$ implica que $y\in P_{s_0s_1\ldots s_n}$. De esta manera $I(x) $ e $I(y)$   coinciden  en los primeros $n+1$ t\'erminos y as\'i, por la proposici\'on \ref{proposicion221}
	\[d[I(x),I(y)]\leq \dfrac{1}{N^n}<\epsilon.\]
\end{proof}
\begin{propo}\label{corrimiento}
	La aplicaci\'on de corrimiento unilateral $\sigma:\sum\rightarrow \sum$  dada por \[\sigma(x_i)_{i\in \mathbb{N}}=(x_{i+1})_{i\in\mathbb{N}}\] es una funci\'on continua.
\end{propo} 

\begin{proof}
	Sea $\epsilon>0$ y sea $s=(s_0s_1s_2\ldots )$, se quiere mostrar que $\sigma $ es continua en $s$.\\
	Puesto que $\epsilon>0$, se elige  $n$ de tal forma que $\dfrac{1}{N^n}<\epsilon.$ Al tomar $\delta=\dfrac{1}{N^{n+1}}$. Si $t$ es un punto en $\sum$ y $d[t,s]<\delta$, entonces por \ref{proposicion221}, se puede tener que para $s_i=t_i$ para $i=0,1,\ldots ,n+1$, es decir, $t=(s_0\ldots s_{n+1}t_{n+2}t_{n+3}\ldots ).$  Luego,  $\sigma(t)=(s_1\ldots s_{n+1}t_{n+2}t_{n+3}\ldots )$ tiene entradas que coinciden con las de $\sigma(s)$ en los primeros $n+1$ lugares. As\'i, se tiene 
	\[d[\sigma(s),\sigma(t)]\leq \dfrac{1}{N^n}<\epsilon\]
	Por lo tanto, $\sigma$ es continua en $s.$
\end{proof}

\begin{defi}\label{semiconjugacion}
	Sea $F$ \[F: X \longrightarrow  X \] una función continua por pedazos y $X$ un espacio topológico y sea \[\sigma: \sum_{F}\longrightarrow \sum_{F},\]
	
	se dice  que $F$ y $\sigma_F$ son semiconjugadas si el diagrama 
	
	\[\xymatrix{ \ar[d]^{h} X \ar[r]^{F}  & X \ar[d]^{h}\\
		\sum_{F}  \ar[r]^{\sigma}      & \sum_{F} 
	}\]
	es conmutativo, es decir, 
	\[\sigma_F\circ h=h\circ F\]
	y la función $h$ es suprayectiva.
\end{defi}

La existencia de una semiconjugación entre dos espacios, busca trasladar información dinámica del espacio de las secuencias permitidas a la dinámica de la aplicación $F: X\rightarrow X$. En el caso de que $h$ sea una biyección se dirá que $F$ es una conjugación con $\sigma$.

\section{Resultados}

\subsection{Espacio de s\'imbolos para $f_\Gamma$}

Para $f_\Gamma$ sobre  $\widehat{\mathbb{R}}$, se considera el alfabeto $\mathcal{A}=\{0,2,4\}$ y el espacio  secuencias de $ \sum_2=\left\{2,0,4\right\}^\mathbb{N}$;  así el itinerario de $x$ bajo $f_\Gamma$ como la sucesi\'on $(I_n(x))\in\sum_2$ como \[I_n(x)=i_0(f^{n-1}_\Gamma(x))\]
donde \[i_0(x)=\left\{\begin{array}{cc}
0, &x\in P^0  \\
2,&x\in P^2\\
4,&x\in P^4.
\end{array}\right.\]

	\begin{figure}[h!]
	\centering
	\begin{center} 
		\includegraphics[scale=0.5]{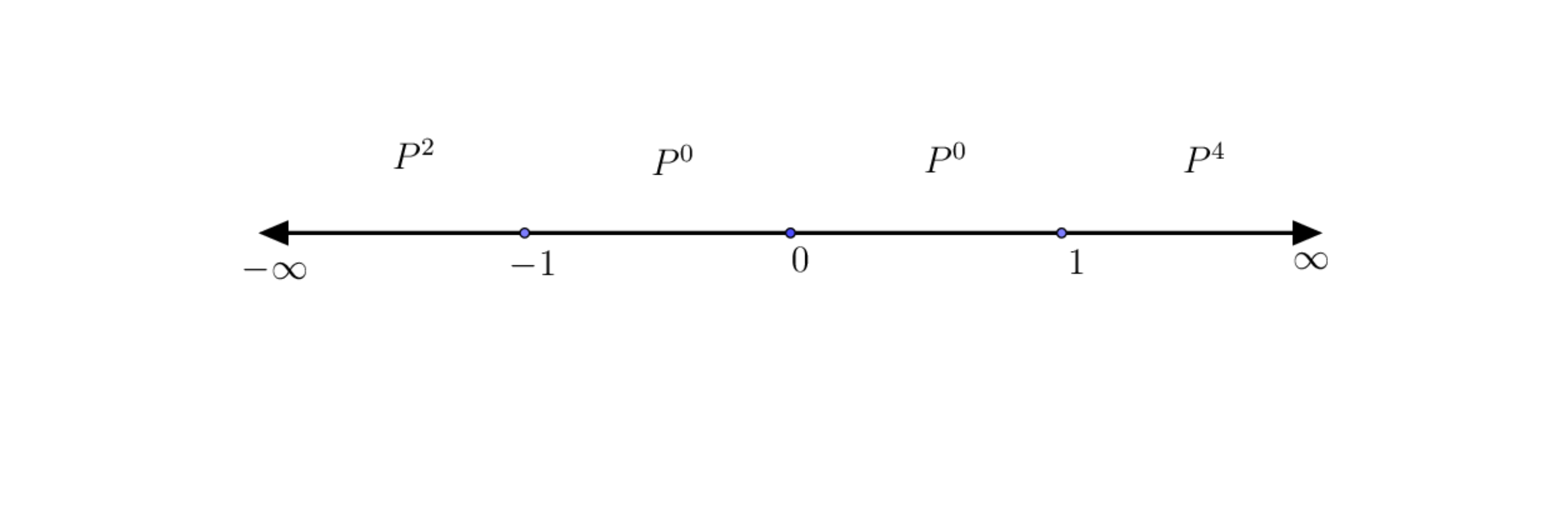}
		\caption{Partici\'on sobre la recta extendida}
		\label{linea}
	\end{center}
\end{figure}

\begin{remarks}
	\begin{enumerate}
		\item La aplicación $f_\Gamma$ tiene $\infty$ como el único punto fijo.
		\item  No existen secuencias que contengan bloques de dos dígitos $24$, $42$ y $00$.
		\item No existen secuencias que contengan bloques de tres dígitos $202$ y $404$.
		\item Se reserva la secuencia $0\cdots 0\cdots $ para  $\infty$.
	\end{enumerate}
\end{remarks}

\begin{defi}
	Se designa por $\sum_{f_{\Gamma}\mid \R}$ a la colección de todas las secuencias permitidas bajo $f_\Gamma$.
\end{defi}

\begin{theorem}
	$f_\Gamma$ es una semiconjugaci\'on con el corrimiento unilateral $\sigma$.
	\[\xymatrix{ \ar[d]^{I}\widehat{\R} \ar[r]^{f_\Gamma}  &  \widehat{\R}\ar[d]^{I}\\
		\sum_{f_\Gamma}  \ar[r]^{\sigma}      & \sum_{f_\Gamma} 
	}\]
\end{theorem}
\begin{proof}
	Para demostrar la semiconjugación según la definición \ref{semiconjugacion}, debemos garantizar que la aplicación $I$ es continua y que el diagrama es conmutativo, por tanto
	\begin{itemize}
		\item En general, las funciones itinerario y de  corrimientos son continuas, esto por las Proposiciones \ref{Itinerario} y \ref{corrimiento} respectivamente.
		\item N\'otese que el diagrama es conmutativo,  ya que:
		
		para $x\in \widehat{\mathbb{R}}$, dicho punto tiene itinerario $s_0s_1s_2\ldots $, por definici\'on 
		\begin{align*}
		x&\in P^{s_0}\\
		f_\Gamma(x)&\in P^{s_1}\\
		f^2_\Gamma(x)&\in P^{s_2}\\
		\vdots &
		\end{align*}
		donde $P^{s_j}$ puede ser  $P^0=(-1,0)\cup (0,1)$, $P^2=(-\infty,-1)$ o $P^4=(1,\infty)$, de donde 
		\[I(f_\Gamma(x))=s_1s_2s_3\ldots \] que es justamente $\sigma(I(x)).$

		\item  Respecto de la sobreyectividad de $I$:
		
		\[\xymatrix{ \ar[d]^{I}\widehat{\R}\ar[r]^{f_\Gamma}  &  \widehat{\R}\ar[d]^{I}\\
			\sum_{f_\Gamma}  \ar[r]^{\sigma}      & \sum_{f_\Gamma} 
		}\]
		Las secuencias en $\displaystyle \sum_{f_\Gamma}$ son periódicas o no periódicas:
		
		Si $(s_i)=\overline{t_0t_1t_2\cdots t_n}$, los s\'imbolos son $2$, $0$ y $4$  y no es posible encontrar las palabras $202$, $404$, $00$, $24$ y $42$; por otro lado, el ciclo periódico debe contener almenos dos ceros.

		Para  $(s_i)=\overline{t_0t_1t_2\cdots t_n}$
		\begin{itemize}
			\item Si $t_0=0$, entonces $t_1=4$ o bien $t_1=2$, por tanto, se puede  asociar la fracción 
			\[ x=-\cfrac{1}{(-1)^sk_1-\cfrac{1}{ (-1)^sk_2\ldots -\frac{1}{ (-1)^sk_l}}}\]
			donde los $k_j$ indican el n\'umero de veces que aparecen los s\'imbolos $2$ o $4$ de manera consecutiva en la secuencia peri\'odica asociada y $s=1$ si $t_{k_j}=2$ o bien $s=2$ si $t_{k_j}=4$.

			\item Si $t_0\neq 0$, entonces $t_0=4$ o bien $t_0=2$, por tanto, es posible asociar la fracción 
			\[ x= (-1)^sk_1-\cfrac{1}{(-1)^sk_2-\cfrac{1}{(-1)^sk_3\ldots -\frac{1}{(-1)^sk_l}}}\]
			donde los $k_j$ indican el n\'umero de veces que aparecen los s\'imbolos $2$ o $4$ de manera consecutiva en la secuencia peri\'odica asociada y $s=1$ si $t_{k_j}=2$ o bien $s=2$ si $t_{k_j}=4$.
			
		\end{itemize}
	\end{itemize}

	Para  $(s_i)=t_0t_1t_2\cdots t_nt_{n+1}\cdots$, sin  patrón periódico.
	
	\begin{itemize}
		\item Los bloques que aparecen son de la forma 
		\[ \underbrace{2\cdots 2}_{m-veces}0\underbrace{4\cdots 4}_{n-veces}\]
		
		se asocia  la fracción 
		\[\cdots -m-\frac{1}{n+\cdots }\]
		
		\item Los bloques que aparecen son de la forma 
		\[ \underbrace{4\cdots 4}_{m-veces}0\underbrace{2\cdots 2}_{n-veces}\]
		
		se asocia  la fracción 
		\[\cdots m-\frac{1}{-n+\cdots }\]
	\end{itemize}
	
\end{proof}

\begin{exam}
	A la secuencia $(s)=\overline{440220}$, se le puede asociar
	el número 
	\[x=2-\cfrac{1}{-2-\cfrac{1}{2+\cdots }},\]  equivalentemente el número como fracción continua
	\[x=2+[-2:2:-2:2:\cdots ].\]
\end{exam}

\begin{exam}
	Por otro lado, la secuencia $(s)=\overline{022044}$
	vendr\'ia del número 
	\[x=-\cfrac{1}{-2-\cfrac{1}{2+\cdots }}\]  equivalentemente el número como fracción continua
	\[x=[0;-2:2:-2:2:-2\cdots ].\]
\end{exam}

\subsection{Estudio de \'orbitas peri\'odicas}
\begin{propo}\label{Pirracionales}
	
	Para $f_\Gamma$ se tiene que: el número $x\in \R$ es peri\'odico bajo $f_\Gamma $ si  y s\'olo si $x=[x_0; x_1:x_2:\ldots ]$ es un irracional cuadrático. 
\end{propo}

\begin{proof} 
	\begin{itemize}
		\item Si $x=0$, por la definici\'on de $f_\Gamma$, $f_z(0)=\infty$ y $f^2_\Gamma(-1)=f^2_\Gamma(1)=\infty$.
		
		\item Si $x\neq 0$ y  $x\in (-1,1)$ es peri\'odico bajo $f_\Gamma $, entonces existe $N\in \mathbb{Z}$ tal que \[f^N_\Gamma(x)=x.\]
		Y puesto que  $f^N_\Gamma$ es la composici\'on de  las transformaciones $S(x)=-\dfrac{1}{x}$ y $T(x)=x+1$, se tiene que existen $M_1,M_2,\cdots, M_n$ de tal forma que 
		
		\[T^{M_1}S\cdots T^{M_{n}}S(x)=f^N_\Gamma(x)=x\] y \[|M_1|+|M_2|+\cdots+|M_n|+n_S=N\]
		donde $n_S$ es el n\'umero de veces en que se aplica la transformación $S(x)=-\cfrac{1}{x}$.

		\begin{align*}
		(T^{M_1}S \ldots T^{M_{n}})\left(-\dfrac{1}{x}\right)&=x\\
		(T^{M_1}S\ldots T^{M_{n-3}}S)\left( -\dfrac{1}{x}+M_{n}\right)&=x.
		\end{align*}
		
		El proceso anterior indica que el n\'umero $x$ tendr\'a la forma
		\[x=-\dfrac{1}{M_1-\dfrac{1}{M_3+\cdots -\dfrac{1}{   \left( -\dfrac{1}{x}+M_{n}\right)       }}}.\]

		Los procesos anteriores indican que el n\'umero $x$ tiene la forma $[M_1,M_2,\cdots M_{n}\cdots]$.
		
		Por otro lado, si $x=[x_0: x_1,x_2,\cdots ]$ es un irracional peri\'odico se tiene 
		\[x=x_0+\dfrac{1}{x_1+\dfrac{1}{x_2+\dfrac{1}{x_3\cdots }}}, \quad x_0\in \mathbb{Z}\]
		Por tanto 
		\begin{align*}
		ST^{-x_0}(x)&=-\left(x_1+\dfrac{1}{x_2+\dfrac{1}{x_3+\cdots }}\right)\\
		T^{-x_1}ST^{-x_0}(x)&=\left(-\dfrac{1}{x_2+\dfrac{1}{x_3+\cdots }} \right)
		\end{align*}

		Siguiendo con el proceso anterior, se tiene que $x$ tiene la siguiente forma 
		$$T^{-x_n}ST^{-x_{n-1}}S\cdots T^{-x_1}ST^{-x_0}=x$$
		lo anterior indica que existe $N\in \mathbb{N} $ de tal forma que $f^N_\Gamma(x)=x$, en este caso \[N=|x_0|+|x_1|+|x_2|+\cdots +|x_n|+n.\]
		
		\item Finalmente, si $x\notin (-1,1)$, existe $M_0\in\mathbb{Z}$ tal que $f^{M_0}_\Gamma(x)\in (-1,1)$ y se aplica el argumento anterior.
	\end{itemize}
\end{proof}

\begin{propo}\label{Pirracionales}
	Si $x=[x_0; x_1:x_2:\ldots ]$ es un irracional no cuadrático entonces el número $x$ tiene una representación en fracción continua infinita no periódica bajo $f_\Gamma$. 
\end{propo}

\begin{proof}
	Sup\'ongase que $x$ es un n\'umero irracional no cuadr\'atico y adem\'as es peri\'odico bajo $f_\Gamma$, entonces existir\'ia $N\in \mathbb{N}$ tal que 
	\begin{align*}
	f^N_\Gamma(x)&=x
	\end{align*}
	es equivalente a resolver  una ecuación de la forma 
	\begin{align*}
	\cfrac{Ax+B}{Cx+D}&=x, \quad AD-BC\neq 0\\
	Ax+B&=(Cx+D)x\\
	Ax+B&=Cx^2+Dx\\
	0&=Cx^2+(D-A)x-B,
	\end{align*}
	
	que conlleva a una ecuación cuadrática, por tanto, otro tipo de irracionales distintos a los irracionales cuadráticos, no pueden ser periódicos bajo $f_\Gamma(x)$.
\end{proof}

\subsection*{Algoritmo de fracciones continuas para $f_\Gamma$}
Fundamentado en las ideas de \cite{GCF}, se presenta a continuación el algoritmo para encontrar la representación en fracción continua de un número real bajo $f_\Gamma$.

\begin{enumerate}
	\item Si $x$ es un entero, $x=[x]$, el algoritmo termina. 
	\item  Si $x$ no es un entero, se siguen los dos pasos siguientes:
		\begin{enumerate}
		\item Sustraemos del número su parte entera y escribimos:
		\[ x=[x]-\cfrac{1}{\frac{-1}{x-[x]}}.\]
		Hacemos $a_0=[x]$ y $r_1=\cfrac{-1}{x-a_0}$.
		\item supongamos que se han completado $k-1$ pasos y obtuvimos los números 
		$a_{k-2}$ y $r_{k-1}$. Encontramos $a_{k-1}$ y $r_k$:
		\[ r_{k-1}=[r_{k-1}]-\cfrac{1}{\frac{-1}{r_{k-1}-[r_{k-1}]}}\]
		Entonces, hacemos $a_{k-1}=[r_{k-1}]$ y $r_k=\cfrac{-1}{r_{k-1}-a_{k-1}}$.
		\item Por tanto, si $x$ es un número racional, este puede escribirse como 
		
		\[\cfrac{p}{q}=a_0-\cfrac{1}{a_1-\cfrac{1}{a_2-\cdots \cfrac{1}{a_{n-2}-\cfrac{1}{a_{n-1}-\cfrac{1}{a_n}}}}} \]
		donde $a_n$ es un entero. Por otro lado,  para n\'umeros irracionales $x$, su representaci\'on como fracci\'on continua es infinita y puede ser peri\'odica o no:
		\[x=a_0-\cfrac{1}{a_1-\cfrac{1}{a_2-\cdots \cfrac{1}{a_{n-2}-\cfrac{1}{a_{n-1}-\cfrac{1}{a_n-\cdots}}}}}\]
	\end{enumerate}
\end{enumerate}

\begin{exam} El número racional $x=\frac{9}{7}$, puede escribirse  como 
	\begin{align*}
	\cfrac{9}{7}&=1-\cfrac{1}{-3-\cfrac{1}{2}},
	\end{align*}
	y esta descomposición viene del  itinerario
	\[40222044\overline{0}\] en $\sum_{f_\Gamma}$.
\end{exam}

\begin{exam}
	Para $x=\sqrt{2}$, se tiene 
	\begin{align*}
	x&=\sqrt{2}=1.41421356237310\cdots \\
	f_\Gamma(\sqrt{2})&=-1+\sqrt{2}=0.41421356237310\cdots\\
	f^2_\Gamma(\sqrt{2})&=-\cfrac{1}{-1+\sqrt{2}}=-(\sqrt{2}+1)=-2.41421356237309\cdots\\
	f^3_\Gamma(\sqrt{2})&=1-\cfrac{1}{-1+\sqrt{2}}=-\sqrt{2}=-1.41421356237310\cdots\\
	f^4_\Gamma(\sqrt{2})&=2-\cfrac{1}{-1+\sqrt{2}}= -\sqrt{2}+1=-0.41421356237310\cdots\\
	f^5_\Gamma(\sqrt{2})&=-\cfrac{1}{2-\cfrac{1}{-1+\sqrt{2}}}= -1+\sqrt{2}=2.41421356237310\cdots\\
	f^6_\Gamma(\sqrt{2})&=1-\cfrac{1}{2-\cfrac{1}{-1+\sqrt{2}}}= \sqrt{2}=1.41421356237310\cdots
	\end{align*}
	y esta descomposici\'on viene del itinerario \[\overline{402204}\] en $\sum_{f_\Gamma}$.
	\begin{figure}[h!]
		\centering
		\includegraphics[scale=.6]{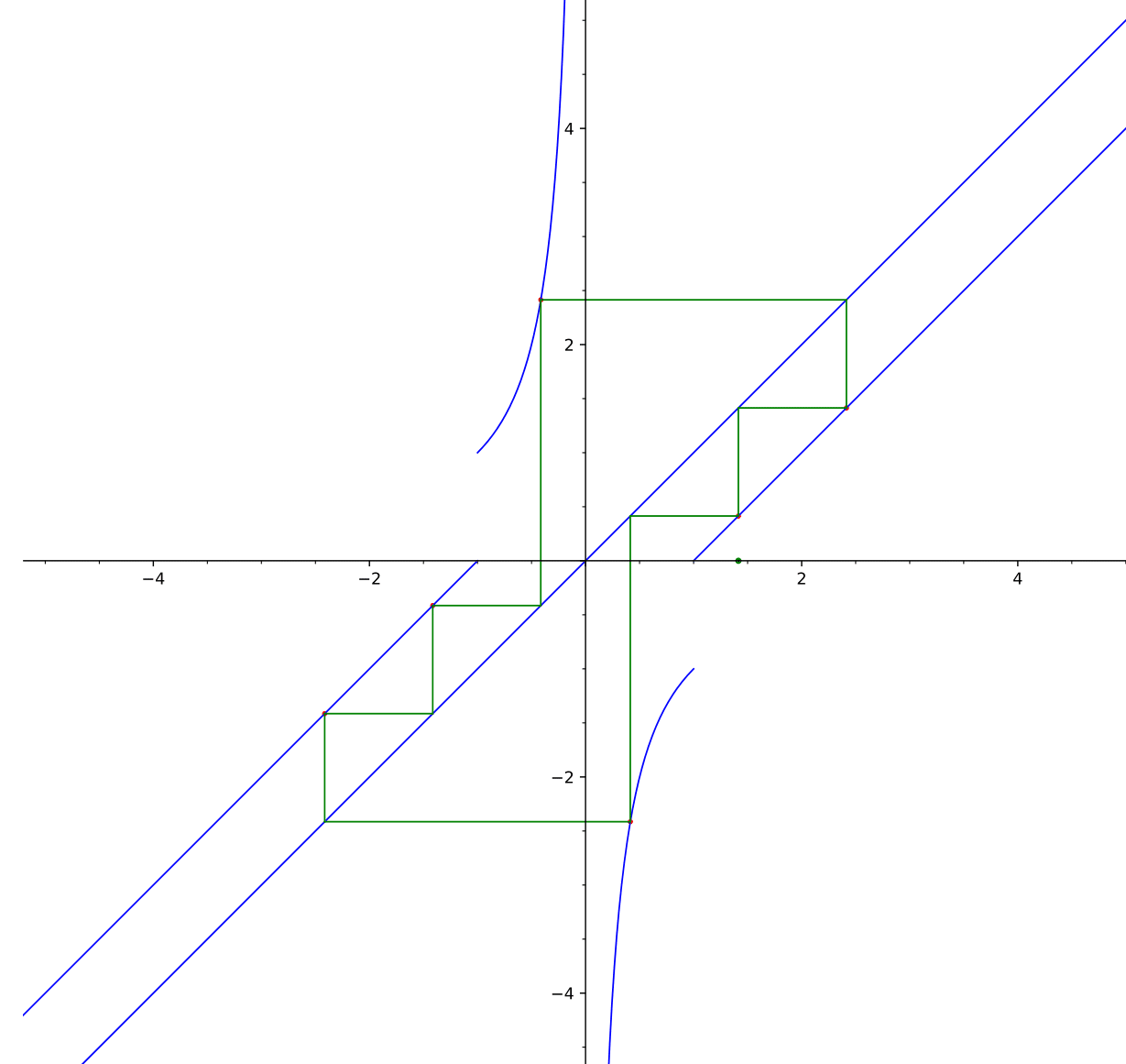}
		\label{fig:my_label}
		\hspace{-5cm}\caption{ Trayectoria de la \'orbita de $x=\sqrt{2}$.}
	\end{figure}
\end{exam}

\section{Estudio numérico de \'orbitas}

\subsection{Trayectorias de números racionales}
\begin{exam}
	Para $x=\frac{4}{3}$, y $x=\frac{5}{4}$ se muestran sus trayectorias 
	
	\begin{figure}[h!]
		\centering
		\includegraphics[width=8cm]{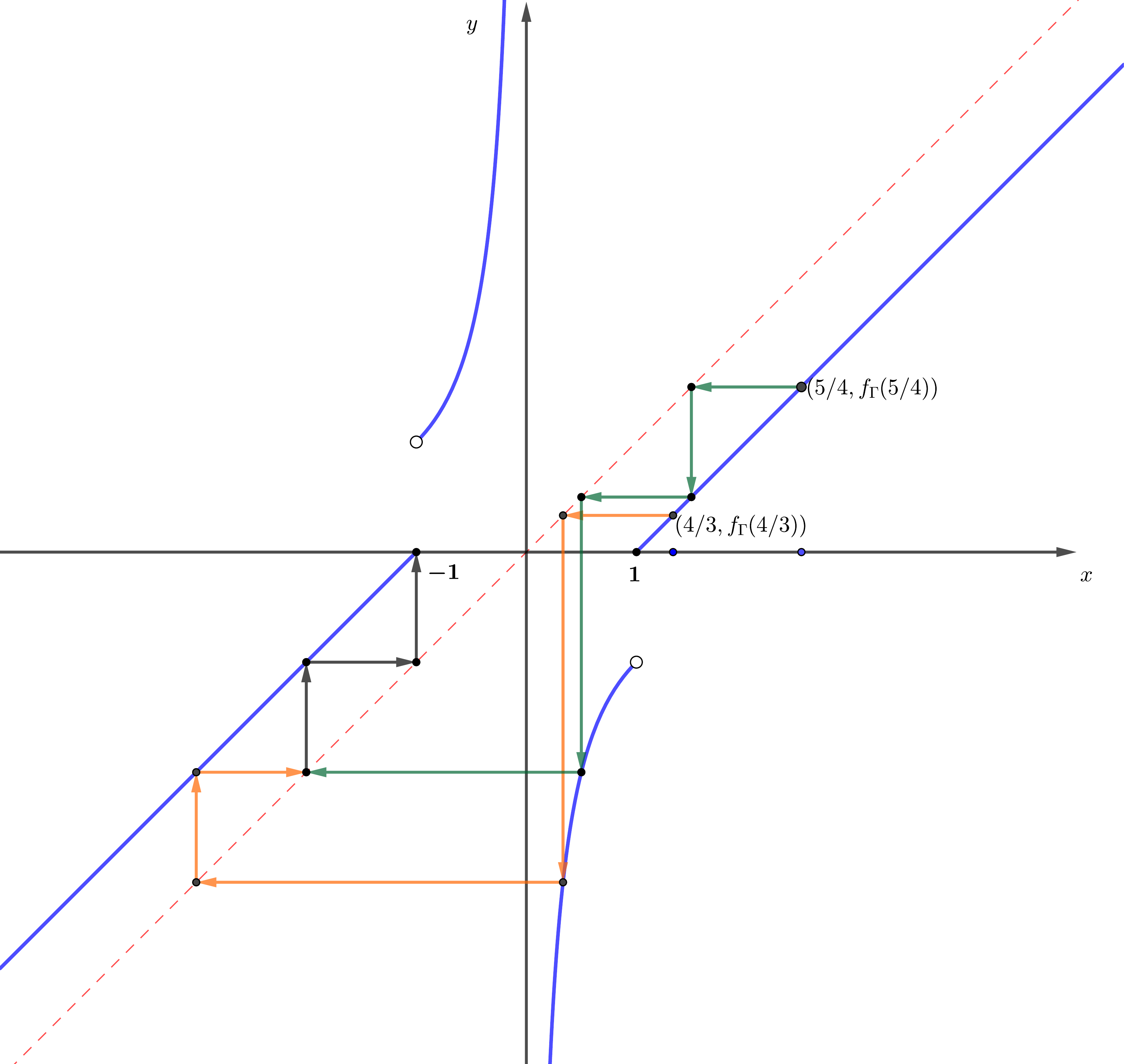}
		\label{racionales}
		\caption{ Trayectorias de la \'orbita de $x=\frac{4}{3}$ y $x=\frac{5}{4}$.}
	\end{figure}
\end{exam}

\subsection{Trayectorias de números irracionales cuadráticos}

La proposici\'on \ref{Pirracionales} garantiza que cualquier n\'umero irracional cuadr\'atico es peri\'odico bajo $f_\Gamma$.

\begin{exam}
	Para $x=\sqrt{2}$ y $x=\cfrac{1+\sqrt{5}}{2}$ se muestran sus trayectorias 
	
	\begin{figure}[h!]
		\centering
		\includegraphics[width=10cm]{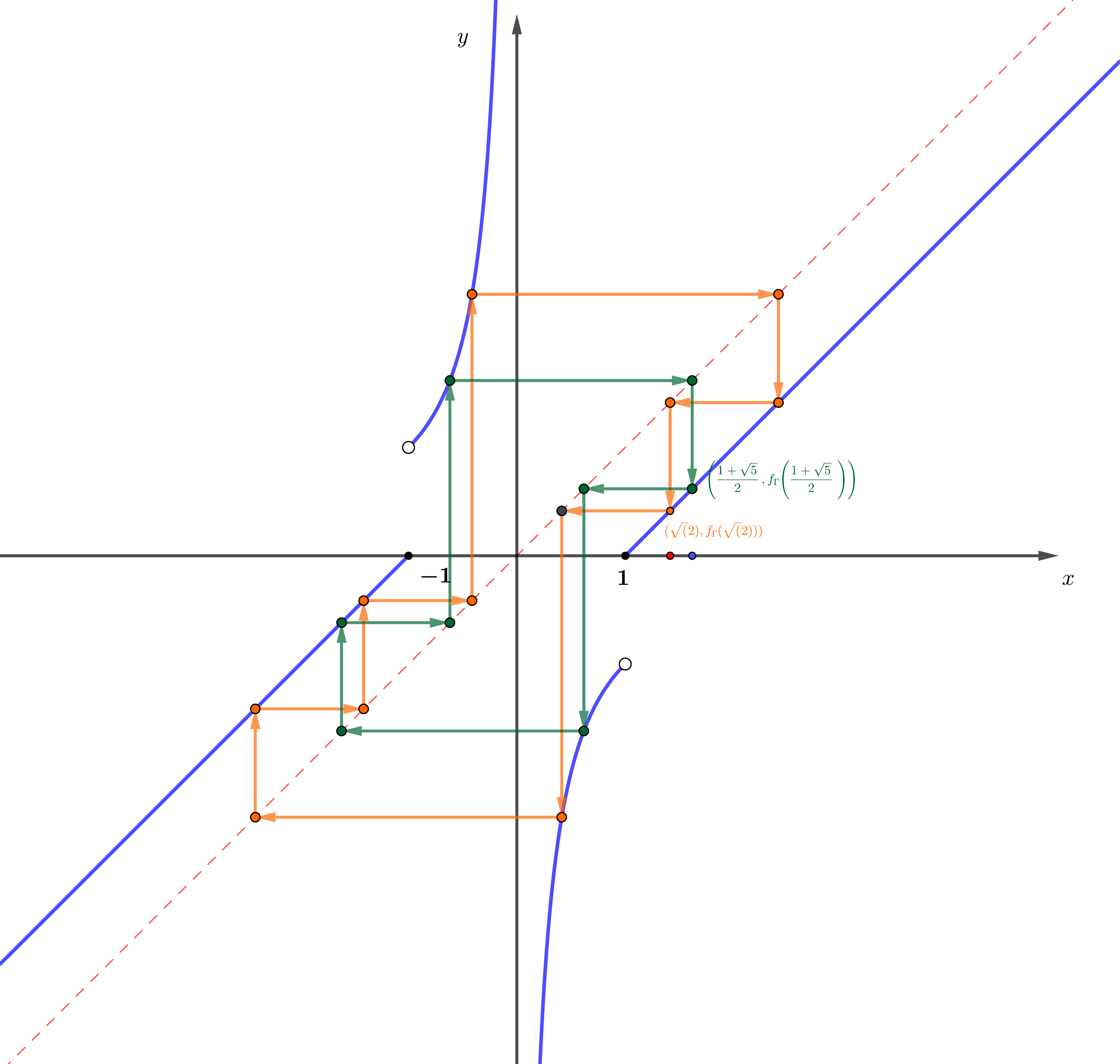}
		\label{irrcuadraticos}
		\caption{ Trayectorias de la \'orbita de $x=\sqrt{2}$ y $x=\left(\frac{1+\sqrt{5}}{2}\right)$.}
	\end{figure}
\end{exam}

\subsection{Trayectorias de números irracionales cúbicos}

\begin{exam}
	Para $\sqrt[3]{2}$, se tiene \\
\begin{multicols}{2}
	\begin{align*}
	x&=\sqrt[3]{2}=1.25992104989487\\
	f_\Gamma(\sqrt[3]{2})&=0.259921049894873\\
	&=\vdots \\
	f^5_\Gamma(\sqrt[3]{2})&=-0.847322101863072\\
	&=\vdots \\
	f^{17}_\Gamma(\sqrt[3]{2})&=-0.22092167902521\\
	&=\vdots\\
	f^{22}_\Gamma(\sqrt[3]{2})&=0.52649103706065
		\end{align*}
			\begin{align*}	&=\vdots \\
	f^{26}_\Gamma(\sqrt[3]{2})&=0.11189244189255\\
	\vdots &=\vdots \\
	f^{35}_\Gamma(\sqrt[3]{2})&=-0.93715413736627\\
	f^{36}_\Gamma(\sqrt[3]{2})&=1.06706032671461\\
	f^{37}_\Gamma(\sqrt[3]{2})&=0.06706032671461\\
	f^{38}_\Gamma(\sqrt[3]{2}) &=-14.9119464367611\\
	\vdots &=\vdots\\
	f^{115}_\Gamma(\sqrt[3]{2})&=-0.0755733032233
	\end{align*}
\end{multicols}
\end{exam}

		\begin{figure}[h!]
			\centering
		\includegraphics[scale=.6]{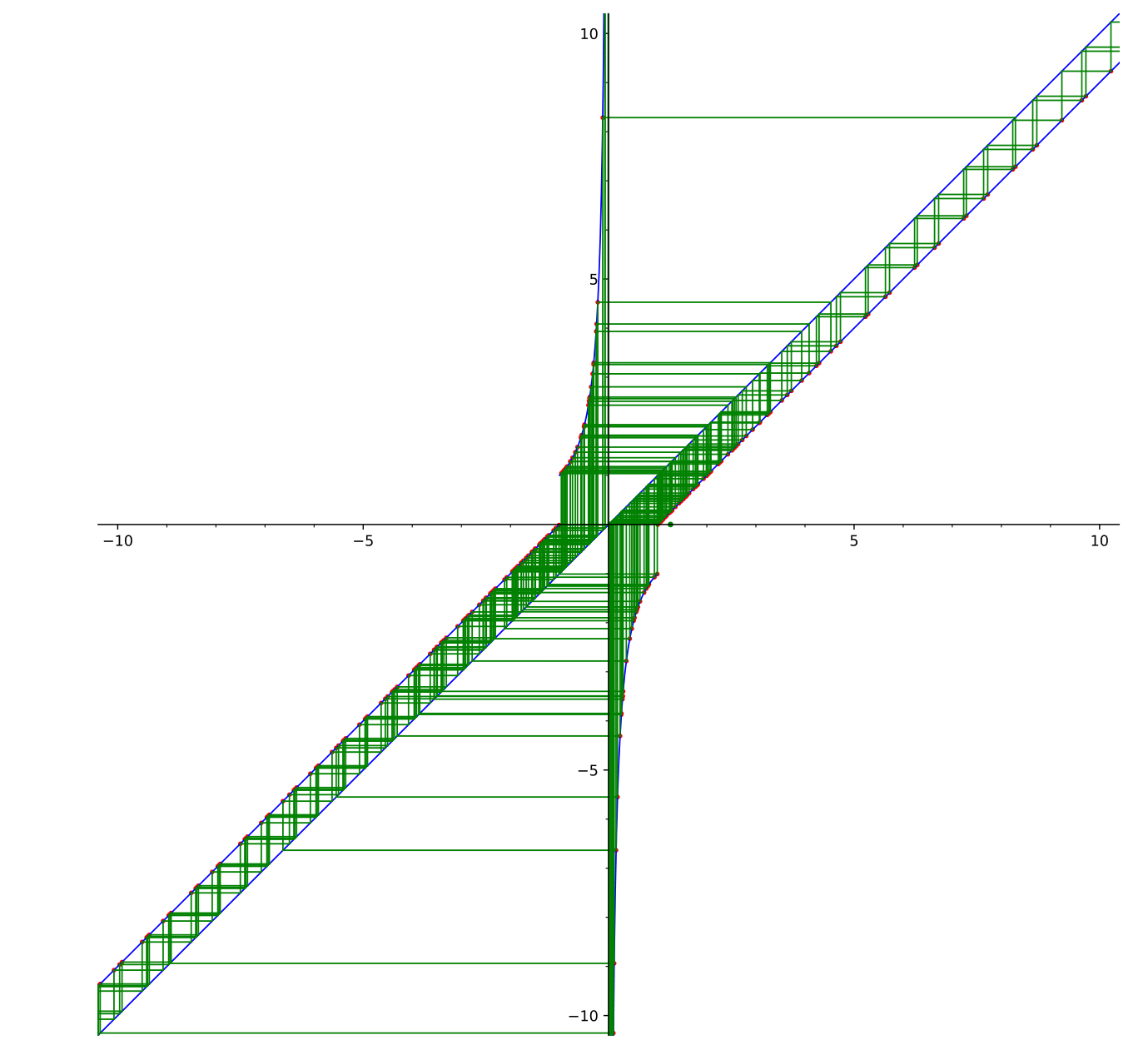}
			\label{raizcubica}
		\caption{ Trayectoria de la \'orbita de $x=\sqrt[3]{2}$.}
	\end{figure}

	El itinerario asociado a $x=\sqrt[3]{2}$ estaría dado por 	
	\[40222204444444444440222220\cdots, \]
	y su representación como fracción continua tiene la forma 
	\[\sqrt[3]{2}=1-\cfrac{1}{-4-\cfrac{1}{12-\cfrac{1}{-5\cdots } }},\]
	la cual es infinita ya que su representaci\'on corresponde a un n\'umero irracional, pero  no puede ser peri\'odica bajo $f_\Gamma$ seg\'un \ref{Pirracionales}.

     \begin{exam}
	Para $\sqrt[3]{3}$, se tiene \\
\begin{multicols}{2}
		\begin{align*}
		x&=\sqrt[3]{3}=1.44224957030741\cdots \\
		f_\Gamma(\sqrt[3]{3})&=0.442249570307408\\
		f^2_\Gamma(\sqrt[3]{3})&=-2.26116669667966\\
		f^3_\Gamma(\sqrt[3]{3})&= -1.26116669667966\\
		f^4_\Gamma(\sqrt[3]{3})&=-0.261166696679657\\
		f^5_\Gamma(\sqrt[3]{3})&= 3.82897211900867
			\end{align*}
				\begin{align*}
		f^6_\Gamma(\sqrt[3]{3})&=2.82897211900867\\
		f^7_\Gamma(\sqrt[3]{3})&=1.82897211900867\\
		f^8_\Gamma(\sqrt[3]{3})&= 0.828972119008669\\
		f^9_\Gamma(\sqrt[3]{3})&= -1.20631318842889\\
		f^{10}_\Gamma(\sqrt[3]{3}) &=-0.20631318842889  \\
		\vdots &=\vdots 
		\end{align*}
	\end{multicols}
	\end{exam}
	
		\begin{figure}[h!]
			\centering
			\includegraphics[scale=.7]{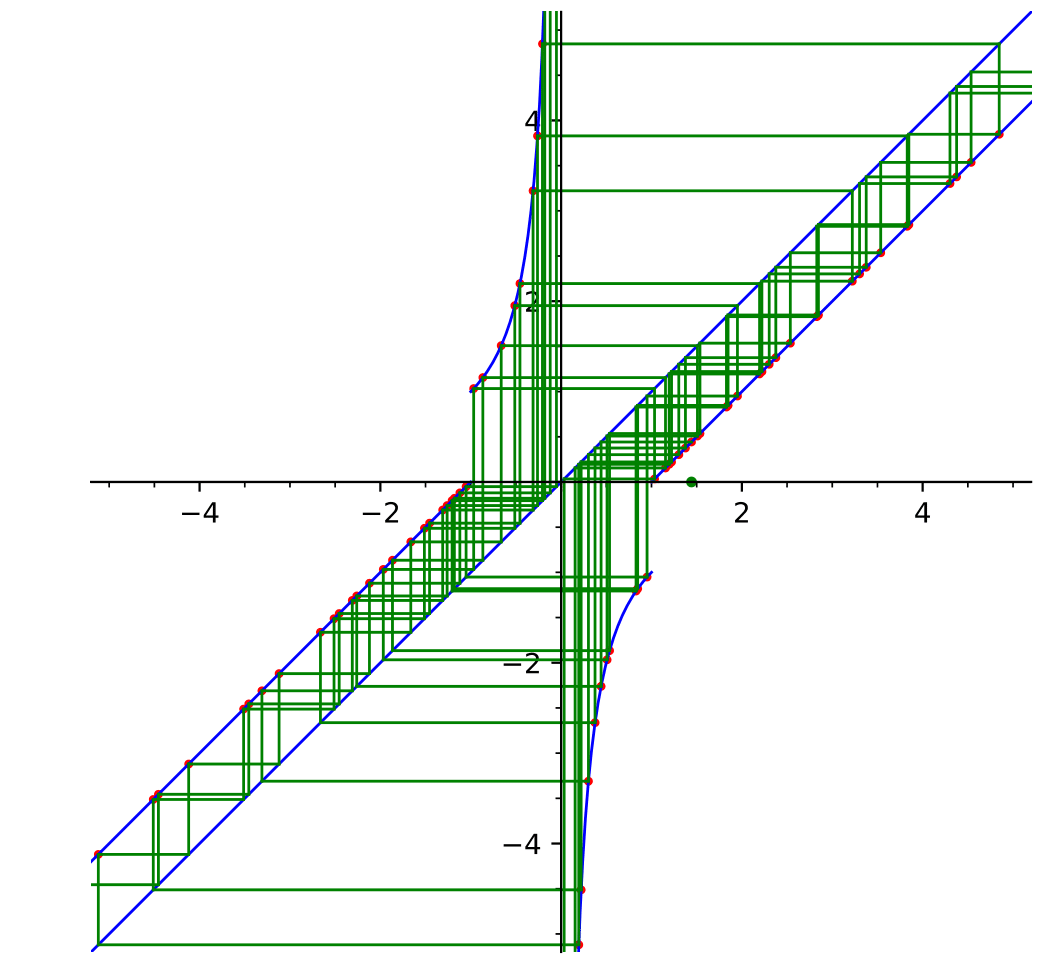}
			\label{cubicatres}
		\caption{ Trayectoria de la \'orbita de $x=\sqrt[3]{3}$.}
	\end{figure}
	El itinerario asociado a $x=\sqrt[3]{3}$ estaría dado por 
	
	\[40220444020\cdots  \]
	y su representación como fracción continua tiene la forma 
	\[\sqrt[3]{3}=1-\cfrac{1}{-2-\cfrac{1}{3-\cfrac{1}{-1\cdots } }},\]
	
	la cual es infinita ya que su representaci\'on corresponde a un n\'umero irracional, pero  no puede ser peri\'odica bajo $f_\Gamma$ seg\'un \ref{Pirracionales}.

\subsection{Las \'orbitas de $\pi$}

\begin{multicols}{2}
	\begin{align*}
	x&=\pi=3.141592653589793\\
	&\vdots \\
	f^3_\Gamma(\pi)&=0.141592653589793\\
	&\vdots\\
	f^{11}_\Gamma(\pi)&=-0.06251330593105209\\
	&\vdots 
	\end{align*}
\begin{align*}
	f^{27}_\Gamma(\pi)&= 0.9965944066841033\\
	\vdots\\
	f^{29}_\Gamma(\pi)&=-0.003417231015000244\\
	f^{30}_\Gamma(\pi)&=292.63459087501246\\
	\vdots &\\
	f^{322}_\Gamma(\pi)&=0.6345908750124636
	\end{align*}
	\end{multicols}
	\begin{figure}[h!]
		\centering
	\includegraphics[scale=.5]{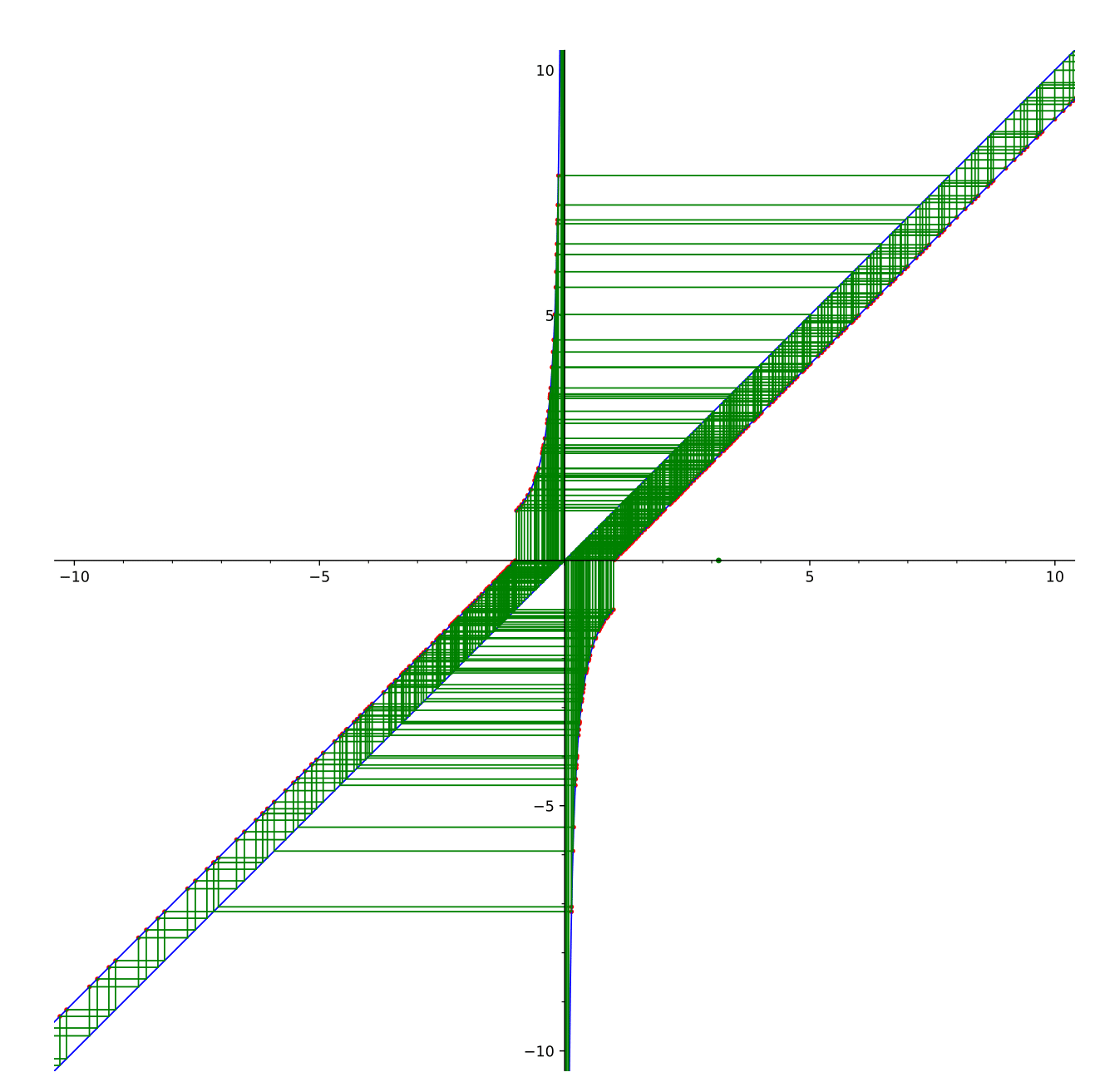}
		\label{orbitapi}
	\caption{ Trayectoria de la \'orbita de $x=\pi$.}
\end{figure}
El itinerario asociado a $x=\pi$ estaría dado por 

\[44402222222044444444444444402\cdots \]
y su representación como fracción continua tiene la forma 
\[\pi=3-\cfrac{1}{-7-\cfrac{1}{15-\cfrac{1}{-1-\cfrac{1}{292\cdots}} }}\]
la cual es infinita ya que su representaci\'on corresponde a un n\'umero irracional, pero  no puede ser peri\'odica bajo $f_\Gamma$ seg\'un \ref{Pirracionales}.

\section{Conclusiones}
\begin{enumerate}
	\item Los n\'umeros reales peri\'odicos bajo $f_\Gamma$ corresponden a los n\'umeros con representaci\'on como fracci\'on continua peri\'odica, es decir,  los  números irracionales cuadráticos, en concordancia con la representaci\'on de n\'umeros irracionales cuadr\'aticos de \cite{GCF}.
	\item Mediante $f_\Gamma$, la representación de un número irracional $x$ no cuadrático, lo cual incluye a los irracionales cúbicos e irracionales trascendentes, tienen representación en fracción continua infinita y por la proposici\'pon \ref{Pirracionales} no pueden ser periódicas.\\ Por otro lado, para cualquier número irracional, es posible asignar un itinerario en $\sum_{\Gamma}$ de manera experimental mediante algoritmos computacionales, es posible aproximar el número irracional $x$ con la fracción continua asociada tal como lo indida el proceso de la sobreyectividad de la Proposición \ref{semiconjugacion}. Por ejemplo, para $x=\sqrt[3]{3}$, su itinerario en $\sum_{f_\Gamma}$ tiene la forma
	
	\[40220444020\cdots  \]
	y su representación como fracción continua tiene la forma 
	\[1.44224957030741\cdots =\sqrt[3]{3}=1-\cfrac{1}{-2-\cfrac{1}{3-\cfrac{1}{-1\cdots } }}=[1;-2:3:-1:\cdots ],\]
	
	al truncar la fracción en el tercer término tendríamos 
	\[1-\cfrac{1}{-2-\cfrac{1}{3-\cfrac{1}{-1 } }}=\cfrac{259}{209}=1.44\cdots, \]
	obteniendo dos cifras de exactitud en el tercer término. Por tanto, pese a no generar una representación periódica como fracción continua, es posible lograr aproximaciones, no solamente para irracionales cúbicos, sino, para cualquier tipo de irracionales.
	
	\item Dada una fracción continua generada por por medio de $f_\Gamma$ de cualquier número irracional $x$, podemos aproximar el valor de un número irracional haciendo cálculos elementales de sumas y divisiones de números enteros. Esto puede ser utilizado de manera didáctica para entender la representación decimal y geométrica de números irracionales en procesos de formación básica y media.
\end{enumerate}

\bibliographystyle{plain}

\end{document}